\documentclass[12pt]{article}

\usepackage{amsmath,amssymb,amsthm,amscd}

\begin{document}

%%%  Abbreviations

\newcommand{\non}{\nonumber}
\newcommand{\scl}{\scriptstyle}
\newcommand{\wt}{\widetilde}
\newcommand{\wh}{\widehat}
\newcommand{\ot}{\otimes}
\newcommand{\fand}{\quad\text{and}\quad}
\newcommand{\Fand}{\qquad\text{and}\qquad}
\newcommand{\ts}{\,}
\newcommand{\tss}{\hspace{1pt}}
\newcommand{\tpr}{t^{\tss\prime}}
\newcommand{\spr}{s^{\tss\prime}}

%%% Roman letters and notation

\newcommand{\U}{ {\rm U}}
\newcommand{\Z}{ {\rm Z}}
\newcommand{\ZY}{ {\rm ZY}}
\newcommand{\Sr}{ {\rm S}}
\newcommand{\Ar}{ {\rm A}}
\newcommand{\Br}{ {\rm B}}
\newcommand{\Cr}{ {\rm C}}
\newcommand{\Ir}{ {\rm I}}
\newcommand{\Jr}{ {\rm J}}
\newcommand{\Prm}{ {\rm P}}
\newcommand{\X}{ {\rm X}}
\newcommand{\Y}{ {\rm Y}}
\newcommand{\SY}{ {\rm SY}}
\newcommand{\Pf}{ {\rm Pf}\ts}
\newcommand{\Hf}{ {\rm Hf}\ts}
\newcommand{\tr}{ {\rm tr}\ts}
\newcommand{\End}{{\rm{End}\ts}}
\newcommand{\Hom}{{\rm{Hom}}}
\newcommand{\sgn}{ {\rm sgn}\ts}
\newcommand{\sign}{ {\rm sign}\ts}
\newcommand{\qdet}{ {\rm qdet}\ts}
\newcommand{\sdet}{ {\rm sdet}\ts}
\newcommand{\Ber}{ {\rm Ber}\ts}
\newcommand{\inv}{ {\rm inv}\ts}
\newcommand{\inva}{ {\rm inv}}
\newcommand{\middle}{ {\rm mid} }
\newcommand{\ev}{ {\rm ev}}
\newcommand{\gr}{ {\rm gr}\ts}

%%% Bb style

\newcommand{\CC}{\mathbb{C}}
\newcommand{\ZZ}{\mathbb{Z}}

%%% Cal style

\newcommand{\Ac}{{\mathcal A}}
\newcommand{\Bc}{{\mathcal B}}
\newcommand{\Oc}{{\mathcal O}}
\newcommand{\Pc}{{\mathcal P}}
\newcommand{\Qc}{{\mathcal Q}}

%%% Gothic style

\newcommand{\Sym}{\mathfrak S}
\newcommand{\h}{\mathfrak h}
\newcommand{\gl}{\mathfrak{gl}}
\newcommand{\oa}{\mathfrak{o}}
\newcommand{\spa}{\mathfrak{sp}}
\newcommand{\osp}{\mathfrak{osp}}
\newcommand{\g}{\mathfrak{g}^{}}
\newcommand{\ggot}{\mathfrak{g}}
\newcommand{\agot}{\mathfrak{a}}
\newcommand{\sll}{\mathfrak{sl}}

%%% Greek letters

\newcommand{\al}{\alpha}
\newcommand{\be}{\beta}
\newcommand{\ga}{\gamma}
\newcommand{\de}{\delta^{}}
\newcommand{\ve}{\varepsilon}
\newcommand{\vp}{\varphi}
\newcommand{\la}{\lambda}
\newcommand{\La}{\Lambda}
\newcommand{\Ga}{\Gamma}
\newcommand{\si}{\sigma}
\newcommand{\ze}{\zeta}
\newcommand{\om}{\omega^{}}

%%% Math bf letters

\newcommand{\rb}{\mathbf r}
\newcommand{\sbf}{\mathbf s}

\renewcommand{\theequation}{\arabic{section}.\arabic{equation}}

\newtheorem{thm}{Theorem}[section]
\newtheorem{lem}[thm]{Lemma}
\newtheorem{prop}[thm]{Proposition}
\newtheorem{cor}[thm]{Corollary}

\theoremstyle{definition}
\newtheorem{defin}[thm]{Definition}
\newtheorem{example}[thm]{Example}

\theoremstyle{remark}
\newtheorem{remark}[thm]{Remark}

\newcommand{\bth}{\begin{thm}}
\renewcommand{\eth}{\end{thm}}
\newcommand{\bpr}{\begin{prop}}
\newcommand{\epr}{\end{prop}}
\newcommand{\ble}{\begin{lem}}
\newcommand{\ele}{\end{lem}}
\newcommand{\bco}{\begin{cor}}
\newcommand{\eco}{\end{cor}}
\newcommand{\bde}{\begin{defin}}
\newcommand{\ede}{\end{defin}}
\newcommand{\bex}{\begin{example}}
\newcommand{\eex}{\end{example}}
\newcommand{\bre}{\begin{remark}}
\newcommand{\ere}{\end{remark}}

\newcommand{\bal}{\begin{aligned}}
\newcommand{\eal}{\end{aligned}}
\newcommand{\beq}{\begin{equation}}
\newcommand{\eeq}{\end{equation}}
\newcommand{\ben}{\begin{equation*}}
\newcommand{\een}{\end{equation*}}

\newcommand{\bpf}{\begin{proof}}
\newcommand{\epf}{\end{proof}}

\def\beql#1{\begin{equation}\label{#1}}

\title{\Large\bf  Verma modules for Yangians}
\author{{Y. Billig, \quad V. Futorny\quad and\quad A. Molev}}

\date{} % 30 April 2005 - start.

\maketitle

\begin{abstract}
We study the Verma modules $M(\mu(u))$ over the Yangian $\Y(\agot)$ associated
with a simple Lie algebra $\agot$. We give necessary and sufficient conditions
for irreducibility of $M(\mu(u))$. Moreover, regarding the simple quotient
$L(\mu(u))$ of $M(\mu(u))$ as an $\agot$-module, we give
necessary and sufficient conditions for finite-dimensionality of the weight
subspaces of $L(\mu(u))$.
\end{abstract}

%%%\vspace{7 mm}
%%%
%%%{\it Key words:}
%%%

\vspace{33 mm}

\noindent
School of Mathematics and Statistics\newline
Carleton University,
1125 Colonel By Drive\newline
Ottawa, Ontario, K1S 5B6, Canada\newline
billig@math.carleton.ca

\vspace{7 mm}

\noindent
Instituto de Matematica e Estatistica\newline
Universidade de S\~ao Paulo\newline
Caixa Postal 66281 CEP 05315-970\newline
S\~ao Paulo, Brazil\newline
futorny@ime.usp.br

\vspace{7 mm}

\noindent
School of Mathematics and Statistics\newline
University of Sydney,
NSW 2006, Australia\newline
alexm@maths.usyd.edu.au

\newpage

\section{Introduction}\label{sec:int}
\setcounter{equation}{0}

For a simple Lie algebra $\agot$ over $\CC$ consider the corresponding
Yangian $\Y(\agot)$; see Drinfeld~\cite{d:ha, d:nr} and the definition
in Section~\ref{sec:def} below. Let $n$ be the rank of the Lie algebra $\agot$.
Given any $n$-tuple $\mu(u)=(\mu_1(u),\dots,\mu_n(u))$
of formal series
\beql{mui}
\mu_i(u)=1+\mu_i^{(0)}u^{-1}+\mu_i^{(1)}u^{-2}+\cdots,
\qquad \mu_i^{(r)}\in\CC,
\eeq
the Verma module $M(\mu(u))$ over $\Y(\agot)$ is defined in
a standard way as the quotient of $\Y(\agot)$ by a left ideal;
see Section~\ref{sec:def}.
The $n$-tuple $\mu(u)$ is called the highest weight of $M(\mu(u))$.
A standard argument shows that $M(\mu(u))$ has a unique simple
quotient $L(\mu(u))$. By a theorem of Drinfeld~\cite{d:nr},
every finite-dimensional irreducible representation of the Yangian
$\Y(\agot)$ is isomorphic to $L(\mu(u))$ for an appropriate
highest weight $\mu(u)$. Moreover, the same theorem gives
necessary and sufficient conditions for the representation $L(\mu(u))$
to be finite-dimensional.

In this paper we are concerned with the necessary and sufficient
conditions on the highest
weight $\mu(u)$ for the Verma module $M(\mu(u))$ to be irreducible.
Our first main result is the following.

\bth\label{thm:redu}
The Verma module $M(\mu(u))$ over the Yangian $\Y(\agot)$ is reducible
if and only if for some index $i\in\{1,\dots,n\}$ the series $\mu_i(u)$
is the Laurent expansion at $u=\infty$ of a rational function in $u$,
\ben
\mu_i(u)=\frac{P(u)}{Q(u)},
\een
where $P(u)$ and $Q(u)$ are monic polynomials in $u$ of the same degree.
\eth

This theorem will be proved in Section~\ref{sec:irr}.
We would like to note a rather unexpected difference between these reducibility
conditions and those for the Verma modules over the Lie algebra $\agot$
itself. It is well-known (see e.g. Dixmier~\cite[Chapter~7]{d:ae})
that the reducibility conditions for a Verma module over $\agot$
involve arbitrary positive roots of $\agot$, while the conditions
of Theorem~\ref{thm:redu} only involve the simple roots
labelled by the indices $i=1,\dots,n$.

Since the Yangian $\Y(\agot)$ contains the universal enveloping algebra
$\U(\agot)$ as a subalgebra, we may regard both $M(\mu(u))$ and
$L(\mu(u))$ as $\agot$-modules. Each of these modules
admits a weight space decomposition with respect to the Cartan subalgebra
of $\agot$. It is immediate from the definition of $M(\mu(u))$ that
all its weight subspaces are infinite-dimensional, except for the one
spanned by the highest vector. The following is our second main theorem.

\bth\label{thm:fin}
All weight subspaces of the $\agot$-module $L(\mu(u))$
are finite-dimensional
if and only if for each index $i\in\{1,\dots,n\}$ the series $\mu_i(u)$
is the Laurent expansion at $u=\infty$ of a rational function in $u$,
\ben
\mu_i(u)=\frac{P_i(u)}{Q_i(u)},
\een
where $P_i(u)$ and $Q_i(u)$ are monic polynomials in $u$ of the same degree.
\eth

We prove this theorem in Section~\ref{sec:fin}.
Note that by Drinfeld's theorem~\cite{d:nr},
the finite-dimensional modules $L(\mu(u))$ correspond to the case
where for each $i=1,\dots,n$ we have $P_i(u)=Q_i(u+d_i)$
for some positive integers $d_i$
defined in the next section;
see also \cite[Theorem~12.1.11]{cp:gq}.
We also note similarity of
Theorem~\ref{thm:fin}
with the recent work of Billig and
Zhao~\cite{bz:wm}
where the finite-dimensionality conditions for the weight subspaces
are found for a wide class of representations of the
exp-polynomial Lie algebras.

\section{Definitions and preliminaries}\label{sec:def}
\setcounter{equation}{0}

As before, we let $\agot$ denote a finite-dimensional simple
Lie algebra over $\CC$ and let $n$ be the rank of $\agot$.
Let $A=(a_{ij})$ be the Cartan matrix of $\agot$. The positive integers
$d_1,\dots,d_n$ are determined by the condition that they are
coprime and the product $DA$ is a symmetric matrix, where
$D=\text{diag}\tss(d_1,\dots,d_n)$. Following Drinfeld~\cite{d:nr}
(see also \cite[Chapter~12]{cp:gq}), we define the {\it Yangian\/}
$\Y(\agot)$ as the associative algebra with generators
$e_i^{(r)}, h_i^{(r)}, f_i^{(r)}$ where $i=1,\dots,n$ and $r=0,1,2,\dots$,
and the following defining relations
\begin{align}
\label{ef}
[h_i^{(r)},h_j^{(s)}]&=0,\qquad\qquad\quad
[e_i^{(r)},f_j^{(s)}]=\delta_{ij}h_i^{(r+s)},\\
[h_i^{(0)},e_j^{(s)}]&=d_i\tss a_{ij}e_j^{(s)},\qquad
[h_i^{(0)},f_j^{(s)}]=-d_i\tss a_{ij}f_j^{(s)},\non\\
[h_i^{(r+1)},e_j^{(s)}]&-[h_i^{(r)},e_j^{(s+1)}]=\frac12 d_ia_{ij}
(h_i^{(r)}e_j^{(s)}+e_j^{(s)}h_i^{(r)}),\non\\
[h_i^{(r+1)},f_j^{(s)}]&-[h_i^{(r)},f_j^{(s+1)}]=-\frac12 d_ia_{ij}
(h_i^{(r)}f_j^{(s)}+f_j^{(s)}h_i^{(r)}),\non\\
[e_i^{(r+1)},e_j^{(s)}]&-[e_i^{(r)},e_j^{(s+1)}]=\frac12 d_ia_{ij}
(e_i^{(r)}e_j^{(s)}+e_j^{(s)}e_i^{(r)}),\non\\
[f_i^{(r+1)},f_j^{(s)}]&-[f_i^{(r)},f_j^{(s+1)}]=-\frac12 d_ia_{ij}
(f_i^{(r)}f_j^{(s)}+f_j^{(s)}f_i^{(r)}),\non\\
\text{\rm Sym}&\ts[e_{i}^{(r_1)},[e_{i}^{(r_2)},\dots,
[e_{i}^{(r_m)},e_{j}^{(s)}]\dots ]]=0,\non\\
\text{\rm Sym}&\ts[f_{i}^{(r_1)},[f_{i}^{(r_2)},\dots,
[f_{i}^{(r_m)},f_{j}^{(s)}]\dots ]]=0,\non
\end{align}
where the last two relations hold for
all pairs $i\ne j$ with $m=1-a_{ij}$, and $\text{\rm Sym}$
denotes symmetrization with respect to the indices $r_1,\dots,r_m$.

The Yangian $\Y(\agot)$ admits a filtration defined by setting the
degree of $e_i^{(r)}, h_i^{(r)}$ and $f_i^{(r)}$ to be equal to $r$.
For any nonnegative integer $r$ we let $\Y(\agot)_r$ denote
the subspace of $\Y(\agot)$ spanned by the monomials of degree at most $r$ in the
generators. We shall denote the associated graded algebra by $\gr \Y(\agot)$.

The universal enveloping algebra $\U(\agot)$ can be identified
with the subalgebra of $\Y(\agot)$ generated by the
elements $e_i^{(0)}, h_i^{(0)}, f_i^{(0)}$ where $i=1,\dots,n$, so that
$\U(\agot)$ coincides with $\Y(\agot)_0$.
Let us set $h_i=h_i^{(0)}$. The linear span $\h$
of the elements $h_1,\dots,h_n$
is a Cartan subalgebra of the Lie algebra
$\agot$. We take the elements $e_i^{(0)}$
to be the simple root vectors of $\agot$
with respect to $\h$
and set $e_{\al_i}=e_i^{(0)}$
and $f_{\al_i}=f_i^{(0)}$ for $i=1,\dots,n$,
where $\al_1,\dots,\al_n$ denote the simple roots. Let $\Delta^+$
denote the set of positive roots.

We shall need an analog of the Poincar\'e--Birkhoff--Witt theorem
for the algebra $\Y(\agot)$ proved by Levendorski\u\i~\cite{l:pb}.
Let $\al=\al_{i_1}+\dots+\al_{i_p}$ be a decomposition of
a positive root $\al\in\Delta^+$ into a sum of simple roots such that
\ben
e_{\al}=[e_{\al_{i_1}},[e_{\al_{i_2}},\dots,
[e_{\al_{i_{p-1}}},e_{\al_{i_p}}]\dots ]]
\een
is a nonzero root vector corresponding to $\al$, and
\ben
f_{\al}=[f_{\al_{i_1}},[f_{\al_{i_2}},\dots,
[f_{\al_{i_{p-1}}},f_{\al_{i_p}}]\dots ]]
\een
is a nonzero root vector corresponding to $-\al$.
Given a decomposition of a nonnegative integer $r$
into a sum of nonnegative integers $r=r_1+\cdots+r_p$ set
\ben
e_{\al}^{(r)}=[e_{i_1}^{(r_1)},[e_{i_2}^{(r_2)},\dots,
[e_{i_{p-1}}^{(r_{p-1})},e_{i_p}^{(r_p)}]\dots ]]
\een
and
\ben
f_{\al}^{(r)}=[f_{i_1}^{(r_1)},[f_{i_2}^{(r_2)},\dots,
[f_{i_{p-1}}^{(r_{p-1})},f_{i_p}^{(r_p)}]\dots ]].
\een
For any positive root $\al$ the images of the elements $e_{\al}^{(r)}$
and $f_{\al}^{(r)}$
in the $r$-th component of the graded algebra $\gr\Y(\agot)$
are independent of the choice of
partition of $r$: if
$\wt e_{\al}^{\ts(r)}$ is an element obtained
using a different partition, then
$e_{\al}^{(r)}-\wt e_{\al}^{\ts(r)}\in\Y(\agot)_{r-1}$.
The same property is shared by the elements $f_{\al}^{(r)}$.
Given any total ordering on the set
\ben
\{e_{\al}^{(r)}\ |\ \al\in\Delta^+,\ r \geqslant 0\}
\cup \{f_{\al}^{(r)}\ |\ \al\in\Delta^+,\ r \geqslant 0\}\cup
\{h_{i}^{(r)}\ |\ i=1,\dots,n,\ r \geqslant 0\},
\een
the ordered monomials in the elements of this set
form a basis of $\Y(\agot)$; see \cite{l:pb}. This implies that
the associated graded algebra $\gr \Y(\agot)$
is isomorphic to the universal enveloping algebra $\U(\agot\ot\CC[x])$.
The images of the elements $e_{\al}^{(r)}$, $f_{\al}^{(r)}$ and
$h_{i}^{(r)}$ in the $r$-th component of $\gr \Y(\agot)$ can be
identified with $e_{\al}x^{r}$, $f_{\al}x^{r}$ and
$h_{i}x^{r}$, respectively.

Given any $n$-tuple $\mu(u)=(\mu_1(u),\dots,\mu_n(u))$
of formal series \eqref{mui} define the {\it Verma module\/}
$M(\mu(u))$ to be the quotient of $\Y(\agot)$ by the left ideal
generated by the elements $e_i^{(r)}$ and $h_i^{(r)}-\mu_i^{(r)}$
for $i=1,\dots,n$ and $r\geqslant 0$. The image $1_{\mu}$
of the element $1\in\Y(\agot)$ in $M(\mu(u))$ is called
the {\it highest vector\/} of $M(\mu(u))$ and the $n$-tuple
$\mu(u)$ is its {\it highest weight\/}. By the Poincar\'e--Birkhoff--Witt theorem
for the algebra $\Y(\agot)$, given any total ordering
on the set
$
\{f_{\al}^{(r)}\ |\ \al\in\Delta^+,\ r \geqslant 0\},
$
the ordered monomials
\beql{ordmon}
f_{\al^{(1)}}^{(r_1)}\cdots f_{\al^{(l)}}^{(r_l)}\ts 1_{\mu},\qquad
l\geqslant 0,\quad r_i\geqslant 0,\quad \al^{(i)}\in\Delta^+
\eeq
form a basis of $M(\mu(u))$. Regarded as an $\agot$-module, $M(\mu(u))$
has a weight space decomposition with respect to $\h$,
\beql{weidec}
M(\mu(u))=\underset{\eta}{\oplus}\ts M(\mu(u))_{\eta},
\eeq
summed over $n$-tuples $\eta=(\eta_1,\dots,\eta_n)$, where
\ben
M(\mu(u))_{\eta}=\{y\in M(\mu(u))\ |\
h_i\ts y=\eta_i\ts y,\qquad i=1,\dots,n\}.
\een
The weight subspace $M(\mu(u))_{\mu^{(0)}}$ with
$\mu^{(0)}=(\mu_1^{(0)},\dots,\mu_n^{(0)})$ is one-dimensional and spanned
by the highest vector $1_{\mu}$. All other nonzero weight subspaces
correspond to the weights $\eta$ of the form
\beql{weta}
\eta=\mu^{(0)}-k_1\al_1-\cdots-k_n\al_n,
\eeq
where the $k_i$ are nonnegative integers, not all of them are zero.
These weight subspaces are infinite-dimensional.

The sum of all submodules of $M(\mu(u))$ which do not contain
the highest vector $1_{\mu}$ is the unique maximal submodule.
We let $L(\mu(u))$ denote the unique irreducible quotient of $M(\mu(u))$.
It inherits the weight space decomposition
\ben
L(\mu(u))=\underset{\eta}{\oplus}\ts L(\mu(u))_{\eta}
\een
with respect to $\h$.

\section{Irreducibility criterion for the Verma module}\label{sec:irr}
\setcounter{equation}{0}

We present the proof of Theorem~\ref{thm:redu} as a sequence
of propositions. We start by considering
the key case of the theorem when $\agot=\sll_2$.
The algebra $\Y(\sll_2)$ has two more presentations in addition to
the one used in the Introduction;
see Drinfeld~\cite{d:nr}. Following \cite{mno:yc},
we shall use the realization of
$\Y(\sll_2)$ as a subalgebra of the Yangian $\Y(\gl_2)$.
The latter is an associative algebra
with generators $t_{ij}^{(1)},\ t_{ij}^{(2)},\dots$ where
$i,j\in\{1,2\}$, and the defining relations
\beql{defyang}
[t^{(r+1)}_{ij}, t^{(s)}_{kl}]-[t^{(r)}_{ij}, t^{(s+1)}_{kl}]=
t^{(r)}_{kj} t^{(s)}_{il}-t^{(s)}_{kj} t^{(r)}_{il},
\eeq
where $r,s=0,1,\dots$ and  $t^{(0)}_{ij}=\delta_{ij}$.
Introducing the formal generating
series
\beql{tiju}
t_{ij}(u) = \delta_{ij} + t^{(1)}_{ij} u^{-1} + t^{(2)}_{ij}u^{-2} +
\dots\in\Y(\gl_2)[[u^{-1}]],
\eeq
we can write \eqref{defyang} in the equivalent form
\beql{defrel}
(u-v)\ts [t_{ij}(u),t_{kl}(v)]=t_{kj}(u)\ts t_{il}(v)-t_{kj}(v)\ts t_{il}(u).
\eeq
The system of relations \eqref{defyang} is also equivalent to the system
\beql{defequiv}
[t^{(r)}_{ij}, t^{(s)}_{kl}] =\sum_{a=1}^{\min(r,s)}
\Big(t^{(a-1)}_{kj} t^{(r+s-a)}_{il}-t^{(r+s-a)}_{kj} t^{(a-1)}_{il}\Big).
\eeq
Observe that the upper summation index $\min(r,s)$ can be replaced
with $r$ or $s$ because the additional sum which may occur is
automatically zero.

Let $\vp(u)$ be any formal power series in $u^{-1}$
with the leading term $1$,
\ben
\vp(u)=1+\vp_1\ts u^{-1}+\vp_2\ts  u^{-2}+\dots,\qquad \vp_i\in\CC.
\een
Then the mapping
\beql{muphi}
t_{ij}(u)\mapsto \vp(u)\ts t_{ij}(u)
\eeq
defines an automorphism of the algebra $\Y(\gl_2)$.
The Yangian $\Y(\sll_2)$ is isomorphic to the subalgebra
of $\Y(\gl_2)$ which consists of the elements
stable under all automorphisms \eqref{muphi}.
We shall identify $\Y(\sll_2)$ with this subalgebra
by the following formulas,
\beql{isom}
\bal
e(u)&= t_{22}(u)^{-1}t_{12}(u),\\
f(u)&= t_{21}(u)\ts t_{22}(u)^{-1},\\
h(u)&= t_{11}(u)\ts t_{22}(u)^{-1}-t_{21}(u)
\ts t_{22}(u)^{-1}t_{12}(u)\ts t_{22}(u)^{-1},
\eal
\eeq
where
\ben
\bal
e(u)&=\sum_{r=0}^{\infty}e^{(r)}_1u^{-r-1},\\
f(u)&=\sum_{r=0}^{\infty}f^{(r)}_1u^{-r-1},\\
h(u)&=1+\sum_{r=0}^{\infty}h^{(r)}_1u^{-r-1}.
\eal
\een
%and we drop the subscript $1$ in the notation of the generators
%of $\Y(\sll_2)$.
We have the following tensor product decomposition
\beql{ygldec}
\Y(\gl_2)=\ZY(\gl_2)\ot  \Y(\sll_2),
\eeq
where $\ZY(\gl_2)$ denotes the center of $\Y(\gl_2)$.
The series
\ben
\bal
\partial\tss(u)&=t_{11}(u)\ts t_{22}(u-1)-t_{21}(u)\ts t_{12}(u-1)\\
  &=1+\partial_1u^{-1}+\partial_2u^{-2}+\dots
\eal
\een
is called the {\it quantum determinant\/}. All its coefficients
$\partial_1,\partial_2,\dots$
are central in $\Y(\gl_2)$, they are
algebraically independent and generate the center
\cite[Theorem~2.13]{mno:yc}.
Note that
the series $h(u)$ can also be written as
\beql{hudu}
h(u)=t_{22}(u)^{-1}t_{22}(u-1)^{-1} \partial\tss(u).
\eeq

The Poincar\'e--Birkhoff--Witt basis
for the algebra $\Y(\gl_2)$ has the following form \cite[Theorem~1.22]{mno:yc}:
given an arbitrary total ordering on the set of
generators
$t^{(r)}_{ij}$, any element of the algebra $\Y(\gl_2)$ can be uniquely written
as a linear combination of the ordered monomials in the generators.

Given any pair of formal series $\la_1(u),\la_2(u)$,
\beql{lai}
\la_i(u)=1+\la_i^{(1)}u^{-1}+\la_i^{(2)}u^{-2}+\cdots,
\qquad \la_i^{(r)}\in\CC,
\eeq
the Verma module $M(\la_1(u),\la_2(u))$ over $\Y(\gl_2)$ is
the quotient of $\Y(\gl_2)$ by the left ideal
generated by the elements $t_{12}^{(r)}$, $t_{11}^{(r)}-\la_1^{(r)}$
and $t_{22}^{(r)}-\la_2^{(r)}$
for $r\geqslant 1$.  By the Poincar\'e--Birkhoff--Witt theorem
for the algebra $\Y(\gl_2)$, a basis of $M(\la_1(u),\la_2(u))$
is formed by the elements
\beql{basisvm}
t_{21}^{(r_1)}\cdots t_{21}^{(r_k)}\ts 1_{\la},\qquad k\geqslant 0,
\qquad 1\leqslant r_1\leqslant\cdots\leqslant r_k,
\eeq
where $1_{\la}$ is the image
of the element $1\in\Y(\gl_2)$ in $M(\la_1(u),\la_2(u))$.
Note that the ordering of the factors $t_{21}^{(r_i)}$ in
\eqref{basisvm} is irrelevant,
because, by the defining relations,
$[t_{21}^{(r)},t_{21}^{(s)}]=0$ for any $r$ and $s$.

\bpr\label{prop:restr}
The restriction of the
$\Y(\gl_2)$-module $M(\la_1(u),\la_2(u))$ to the subalgebra $\Y(\sll_2)$
is isomorphic to the Verma module $M(\mu(u))$, where the highest weight
$\mu(u)$ is given by $\mu(u)=\la_1(u)/\la_2(u)$.
\epr

\bpf
By the formulas \eqref{isom}, we have
\ben
e(u)\ts 1_{\la}=0,\qquad h(u)\ts 1_{\la}=\mu(u)\ts 1_{\la}
\een
for $\mu(u)=\la_1(u)/\la_2(u)$. Therefore, we have a
$\Y(\sll_2)$-homomorphism $\Phi: M(\mu(u))\to M(\la_1(u),\la_2(u))$
such that $\Phi(1_{\mu})=1_{\la}$. In order to see that $\Phi$ is surjective
let us verify that $M(\la_1(u),\la_2(u))$ is generated by $1_{\la}$
as a $\Y(\sll_2)$-module. Indeed, this follows from the decomposition
\eqref{ygldec} because the elements of $\ZY(\gl_2)$ act on $M(\la_1(u),\la_2(u))$
as scalar operators. Furthermore, the Verma module $M(\mu(u))$ over $\Y(\sll_2)$
with $\mu(u)=\la_1(u)/\la_2(u)$ can be extended to a module over
the Yangian $\Y(\gl_2)$ by defining the action of the quantum determinant
$\partial\tss(u)$ on $M(\mu(u))$ to be the scalar multiplication
by the series $\la_1(u)\la_2(u-1)$. Then by \eqref{isom} and \eqref{hudu}
we have
\ben
t_{12}(u)\ts 1_{\mu}=0,\qquad
t_{ii}(u)\ts 1_{\mu}=\la_i(u)\ts 1_{\mu},\quad i=1,2.
\een
Hence, we have a
$\Y(\gl_2)$-homomorphism $\Psi: M(\la_1(u),\la_2(u))\to M(\mu(u))$
such that $\Psi(1_{\la})=1_{\mu}$. The composition map $\Psi\circ\Phi$
is the identity map on $M(\mu(u))$ which shows that $\Phi$ is injective.
\epf

\bco\label{cor:irr}
The Verma module $M(\la_1(u),\la_2(u))$ over $\Y(\gl_2)$ is irreducible
if and only if the Verma module
$M(\mu(u))$ with $\mu(u)=\la_1(u)/\la_2(u)$
over $\Y(\sll_2)$ is irreducible.
\eco

\bpf
This is immediate from Proposition~\ref{prop:restr}.
\epf

\bpr\label{prop:ratio}
Suppose that formal series $\la_1(u)$ and $\la_2(u)$
given by \eqref{lai} are such that the ratio $\la_1(u)/\la_2(u)$
is the Laurent expansion at $u=\infty$ of a rational function in $u$.
Then the Verma module $M(\la_1(u),\la_2(u))$ over $\Y(\gl_2)$ is reducible.
\epr

\bpf
We can write
\ben
\frac{\la_1(u)}{\la_2(u)}=\frac{P(u)}{Q(u)},
\een
where $P(u)$ and $Q(u)$ are polynomials in $u$. Since the Laurent series
$\la_1(u)/\la_2(u)$ does not contain positive powers of $u$ and has
constant term $1$, we may assume that the polynomials
$P(u)$ and $Q(u)$ are monic of the same degree,
say, $p$. Note that the composition of the action
of $\Y(\gl_2)$ on $M(\la_1(u),\la_2(u))$ with an automorphism
of the form \eqref{muphi} defines a representation of
$\Y(\gl_2)$ on the vector space $M(\la_1(u),\la_2(u))$
which is isomorphic to the Verma module
$M(\vp(u)\la_1(u),\vp(u)\la_2(u))$. Obviously, this Verma module
is reducible or irreducible simultaneously with $M(\la_1(u),\la_2(u))$.
Hence, taking
$
\vp(u)=\la_2(u)^{-1}\ts u^{-p}\tss Q(u),
$
we may assume without loss of generality that both $\la_1(u)$ and $\la_2(u)$
are polynomials in $u^{-1}$ of degree $\leqslant p$.
Consider the vector subspace $K$ of $M(\la_1(u),\la_2(u))$ which is spanned
by the vectors of the form \eqref{basisvm} where at least one of the
indices $r_i$ exceeds $p$. We claim that $K$ is a submodule of
$M(\la_1(u),\la_2(u))$. Indeed, it is obvious that $K$ is stable under
the action of the elements $t_{21}^{(r)}$. Let us now verify
by induction on $k\geqslant 1$ that for any $r\geqslant 1$ and
any positive integers $r_i$ with $r_k\geqslant p+1$ the element
\beql{toneone}
t_{11}^{(r)}\tss t_{21}^{(r_1)}\cdots t_{21}^{(r_k)}\ts 1_{\la}
\eeq
lies in $K$. By the defining relations
\eqref{defequiv}, we have
\beql{defoneone}
[t^{(r)}_{11}, t^{(r_1)}_{21}] =\sum_{a=1}^{\min(r,r_1)}
\Big(t^{(a-1)}_{21} t^{(r+r_1-a)}_{11}-t^{(r+r_1-a)}_{21}
t^{(a-1)}_{11}\Big).
\eeq
If $k=1$ then $r_1\geqslant p+1$ and so $r+r_1-a\geqslant p+1$.
Since $t_{11}^{(r+r_1-a)}\ts 1_{\la}=0$, the statement
is true. Suppose now that $k\geqslant 2$.
The statement now follows from \eqref{defoneone} and
the induction hypothesis.

The same argument shows that the subspace $K$
is stable under the action of the elements
$t^{(r)}_{22}$ with $r\geqslant 1$. Here we use the following consequence
of the relations \eqref{defequiv} instead of \eqref{defoneone}:
\beql{deftwotwo}
[t^{(r)}_{22}, t^{(r_1)}_{21}] =\sum_{a=1}^{\min(r,r_1)}
\Big(t^{(r+r_1-a)}_{21}t^{(a-1)}_{22}-
t^{(a-1)}_{21} t^{(r+r_1-a)}_{22}\Big).
\eeq

Finally, let us verify that the subspace $K$
is stable under the action of the elements
$t^{(r)}_{12}$ with $r\geqslant 1$. We shall show by induction on
$k\geqslant 1$ that all elements of the form
\beql{tonetwo}
t_{12}^{(r)}\tss t_{21}^{(r_1)}\cdots t_{21}^{(r_k)}\ts 1_{\la},
\eeq
where the $r_i$ are positive integers and
$r_k\geqslant p+1$, lie in $K$.
By \eqref{defequiv} we have
\beql{defonetwo}
[t^{(r)}_{12}, t^{(r_1)}_{21}] =\sum_{a=1}^{\min(r,r_1)}
\Big(t^{(a-1)}_{22} t^{(r+r_1-a)}_{11}-t^{(r+r_1-a)}_{22}
t^{(a-1)}_{11}\Big).
\eeq
Hence, if $k=1$ the statement is true because $t^{(r)}_{12}\ts 1_{\la}=0$
and $t^{(s)}_{11}\ts 1_{\la}=t^{(s)}_{22}\ts 1_{\la}=0$ for $s\geqslant p+1$.
Suppose now that $k\geqslant 2$.
Due to \eqref{defonetwo},
the statement follows from the induction hypothesis and
the fact that $K$ is stable under the action of the elements
$t^{(r)}_{11}$ and $t^{(r)}_{22}$ with $r\geqslant 1$.
\epf

\ble\label{lem:recrel}
Let
\beql{nuu}
\nu(u)=1+\nu^{(1)}u^{-1}+\nu^{(2)}u^{-2}+\cdots
\eeq
be a formal series in $u^{-1}$ with complex coefficients.
Suppose that there exist a positive integer $N$ and
complex numbers $c_0,\dots,c_m$, not all zero, such that
the coefficients $\nu^{(r)}$ of the formal series
satisfy the recurrence relation
\beql{recrel}
c_0\ts\nu^{(r)}+c_1\ts\nu^{(r+1)}+\cdots+c_m\ts\nu^{(r+m)}=0,
\eeq
for all $r\geqslant N$. Then $\nu(u)$
is the Laurent expansion at $u=\infty$ of a rational function
$P(u)/Q(u)$ in $u$,
where $P(u)$ and $Q(u)$ are monic polynomials in $u$
of the same degree.
\ele

\bpf
Let us set
\ben
\wt \nu(u)=\nu^{(N)}+\nu^{(N+1)}u^{-1}+\cdots.
\een
Then \eqref{recrel} implies
\ben
\wt \nu(u)\ts(c_0+c_1\ts u+\cdots+c_m\ts u^m)=b_1\ts u+\cdots+b_m\ts u^m
\een
for some coefficients $b_i$. Hence the series
\ben
\nu(u)=1+\nu^{(1)}u^{-1}+\nu^{(2)}u^{-2}+\cdots+\nu^{(N-1)}u^{-N+1}
+u^{-N}\ts \wt \nu(u)
\een
clearly has the desired form.
\epf

\bpr\label{prop:reducib}
Suppose that
the Verma module $M(\la_1(u),\la_2(u))$ over $\Y(\gl_2)$ is reducible.
Then the ratio $\la_1(u)/\la_2(u)$
is the Laurent expansion at $u=\infty$ of a rational function in $u$.
\epr

\bpf
By twisting the action of $\Y(\gl_2)$ on $M(\la_1(u),\la_2(u))$
by the automorphism \eqref{muphi} with $\vp(u)=\la_1(u)^{-1}$,
we obtain a module over $\Y(\gl_2)$ which is isomorphic to
the Verma module $M(1,\nu(u))$ with $\nu(u)=\la_2(u)/\la_1(u)$.
We shall be proving that if the Verma module
$M(1,\nu(u))$ is reducible then
the series $\nu(u)$ is the Laurent expansion
at $u=\infty$ of a rational function in $u$.

Clearly, $M(1,\nu(u))$ is the direct sum
of its $\sll_2$-weight subspaces,
\ben
M(1,\nu(u))=\underset{\eta}{\oplus}\ts M(1,\nu(u))_{\eta},
\een
where
\ben
M(1,\nu(u))_{\eta}=\{y\in M(1,\nu(u))\ |\
(t_{11}^{(1)}-t_{22}^{(1)})\ts y=\eta\ts y\}.
\een
The nonzero weight subspaces correspond to the weights $\eta=-\nu^{(1)}-2k$,
where $k$ is a nonnegative integer and the coefficients
of the series $\nu(u)$ are defined as in \eqref{nuu}.
The monomials \eqref{basisvm} with fixed $k$ form a basis of
$M(1,\nu(u))_{\eta}$ with $\eta=-\nu^{(1)}-2k$.
If $K$ is a nontrivial submodule of $M(1,\nu(u))$ then $K$ inherits
the weight space decomposition,
\ben
K=\underset{\eta}{\oplus}\ts K_{\eta},\qquad
K_{\eta}=K\cap M(1,\nu(u))_{\eta}.
\een
Take the minimum positive integer $k$ such that $K_{\eta}\ne 0$
for $\eta=-\nu^{(1)}-2k$. Then any nonzero vector $\ze\in K_{\eta}$
has the property
$
t_{12}(u)\ts \ze=0
$
since otherwise $t_{12}^{(r)}\ts\ze$ would be a nonzero element of
the subspace $K_{\eta+2}$ for certain $r\geqslant 1$.
Write
\beql{ze}
\ze=\sum_{\rb}c_{\rb}\ts t_{21}^{(r_1)}\cdots t_{21}^{(r_k)}\ts 1_{\la},
\eeq
summed over a finite set of
$k$-tuples $\rb=(r_1,\dots,r_k)$
with $1\leqslant r_1\leqslant\cdots\leqslant r_k$.
Calculating $t_{12}^{(r)}\ts\ze$ as a linear
combination of the basis vectors of $M(1,\nu(u))$ we get a family of
linear relations on the coefficients $\nu^{(s)}$. It will be sufficient
to demonstrate that at least one of these relations has the form
\eqref{recrel} where not all coefficients are zero. The application
of Lemma~\ref{lem:recrel} will then complete the argument.

The action of $t_{12}^{(r)}$ on a basis monomial is calculated by
\beql{tonetwor}
\bal
t_{12}^{(r)}\ts t_{21}^{(r_1)}\cdots t_{21}^{(r_k)}\ts 1_{\la}
{}&={}\sum_{i=1}^k t_{21}^{(r_1)}\cdots t_{21}^{(r_{i-1})}\\
{}&\times{}\sum_{a_i=1}^{r_i}\Big(t^{(a_i-1)}_{22}
t^{(r+r_i-a_i)}_{11}-t^{(r+r_i-a_i)}_{22}
t^{(a_i-1)}_{11}\Big)\ts t_{21}^{(r_{i+1})}\cdots
t_{21}^{(r_k)}\ts 1_{\la}.
\eal
\eeq
In order to write this expression as a linear combination
of basis monomials, we proceed by induction
on $k$ with the use of similar
formulas for the actions of $t_{11}^{(r)}$ and $t_{22}^{(r)}$ with
$r\geqslant 1$
on the basis monomials,
\ben
\bal
t_{11}^{(r)}\ts t_{21}^{(r_1)}\cdots t_{21}^{(r_k)}\ts 1_{\la}
{}&={}\sum_{i=1}^k t_{21}^{(r_1)}\cdots t_{21}^{(r_{i-1})}\\
{}&\times{}\sum_{a_i=1}^{r_i}\Big(t^{(a_i-1)}_{21}
t^{(r+r_i-a_i)}_{11}-t^{(r+r_i-a_i)}_{21}
t^{(a_i-1)}_{11}\Big)\ts t_{21}^{(r_{i+1})}\cdots
t_{21}^{(r_k)}\ts 1_{\la}
\eal
\een
and
\ben
\bal
t_{22}^{(r)}\ts t_{21}^{(r_1)}\cdots t_{21}^{(r_k)}\ts 1_{\la}
{}&={}\nu^{(r)}\ts t_{21}^{(r_1)}\cdots t_{21}^{(r_k)}\ts 1_{\la}+
\sum_{i=1}^k t_{21}^{(r_1)}\cdots t_{21}^{(r_{i-1})}\\
{}&\times{}\sum_{a_i=1}^{r_i}\Big(t^{(r+r_i-a_i)}_{21}
t^{(a_i-1)}_{22}-t^{(a_i-1)}_{21}
t^{(r+r_i-a_i)}_{22}\Big)\ts t_{21}^{(r_{i+1})}\cdots
t_{21}^{(r_k)}\ts 1_{\la},
\eal
\een
where we have used our assumptions on the highest weight, so that
\ben
t_{11}^{(r)}\ts 1_{\la}=0\Fand t_{22}^{(r)}\ts 1_{\la}=
\nu^{(r)}\ts 1_{\la},\qquad
r\geqslant 1.
\een
Using the above formulas, we write \eqref{tonetwor} as a linear
combination of the monomials
$t_{21}^{(s_1)}\cdots t_{21}^{(s_{k-1})}\ts 1_{\la}$.
Suppose that $N$ is the maximum sum of the indices $r_1+\cdots+r_k$
such that the corresponding monomial occurs in \eqref{ze}
with a nonzero coefficient $c_{\rb}$. It will be sufficient
for our purposes to consider only values $r\geqslant N$
and take the coefficient in the expansion
of $t_{12}^{(r)}\ts \ze$ at a monomial of the form
$t_{21}^{(s_1)}\cdots t_{21}^{(s_{k-1})}\ts 1_{\la}$
with $1\leqslant s_1\leqslant\cdots\leqslant s_{k-1}$ and
satisfying $s_1+\cdots+s_{k-1}\leqslant N-1$.
It is not difficult to see
from the above formulas that this coefficient can be written as
a linear combination of $\nu^{(r)},\nu^{(r+1)},\dots,\nu^{(r+N-k)}$.
Moreover, the coefficient at each $\nu^{(r+i)}$
is, in its turn, a linear combination of
the coefficients $c_{\rb}$, and this combination is independent of
the value of the index $r\geqslant N$.
In order to apply Lemma~\ref{lem:recrel}, we need to verify that
at least one of the coefficients at the $\nu^{(r+i)}$ is nonzero.
Let us calculate
the coefficient at $\nu^{(r+s-1)}$ where $s$ is such that
$s_1+\cdots+s_{k-1}+s=N$. This coefficient can only arise from
the expansion of \eqref{tonetwor} with $r_1+\cdots+r_k=N$.
Furthermore, for each $i=1,\dots,k$ the value of the summation
index $a_i$ in that formula must be equal to $1$. Hence, keeping only the terms
contributing to the desired coefficient, we can rewrite \eqref{tonetwor}
as
\ben
t_{12}^{(r)}\ts t_{21}^{(r_1)}\cdots t_{21}^{(r_k)}\ts 1_{\la}
{}\equiv{}-\sum_{i=1}^k \nu^{(r+r_i-1)}
t_{21}^{(r_1)}\cdots t_{21}^{(r_{i-1})}\ts t_{21}^{(r_{i+1})}\cdots
t_{21}^{(r_k)}\ts 1_{\la}.
\een
Thus, the coefficient in question is
\beql{coeffcs}
-c_{\sbf^{(1)}}-c_{\sbf^{(2)}}-\cdots-c_{\sbf^{(k)}},
\eeq
where $\sbf^{(i)}=(s_1,\dots,s_{i-1},s,s_{i},\dots,s_{k-1})$
with $s_{i-1}\leqslant s\leqslant s_{i}$.
We may write \eqref{coeffcs} as $-(m+1)\ts c_{\sbf^{(i)}}$,
where $m$ is the number of indices $s_j$ equal to $s$.
However, by our assumption, there exists a $k$-tuple $\rb=(r_1,\dots,r_k)$
with $r_1+\cdots+r_k=N$ and $c_{\rb}\ne 0$.
So, choosing the parameters $s$ and $s_1,\dots,s_{k-1}$ in such a way
that $c_{\sbf^{(i)}}\ne 0$ for some $i$, we complete the proof.
\epf

Corollary~\ref{cor:irr} and Propositions~\ref{prop:ratio} and
\ref{prop:reducib} imply the following corollary which is
a particular case of Theorem~\ref{thm:redu} for $\agot=\sll_2$.

\bco\label{cor:sltwo}
The Verma module $M(\mu(u))$ over the Yangian $\Y(\sll_2)$ is reducible
if and only if the series $\mu(u)$
is the Laurent expansion at $u=\infty$ of a rational function in $u$.
\qed
\eco

\bco\label{cor:weaksing}
Let $M(\mu(u))$ be the Verma module over the Yangian $\Y(\sll_2)$,
where $\mu(u)=P(u)/Q(u)$ for monic polynomials
$P(u)$ and $Q(u)$ in $u$ of degree $p$.
Then for any $s\geqslant p$ there exist constants
$c_0,\dots,c_{s-1}$ such that the vector
\ben
\ze=c_0\ts f_{1}^{(0)}\tss 1_{\mu}+\cdots+c_{s-1}\ts f_{1}^{(s-1)}\tss 1_{\mu}
+f_{1}^{(s)}\tss 1_{\mu},
\een
satisfies $e(u)\ts \ze=0$.
\eco

\bpf
By Proposition~\ref{prop:restr}, the $\Y(\sll_2)$-module $M(\mu(u))$ is
isomorphic to the restriction of a $\Y(\gl_2)$-module
$M(\la_1(u),\la_2(u))$ such that $\la_1(u)/\la_2(u)=\mu(u)$.
As in the proof of
Proposition~\ref{prop:ratio}, we can find a series $\vp(u)$
such that $\vp(u)\la_1(u)$ and $\vp(u)\la_2(u)$
are polynomials in $u^{-1}$ of degree $\leqslant p$.
For any $s\geqslant p$
the vector $\ze=t_{21}^{(s+1)}\ts 1_{\la}$ belongs to the submodule $K$
introduced in that proof. Hence, $t_{12}^{(r)}\ts\ze=0$ for all
$r\geqslant 1$, since otherwise the highest vector of
the Verma module would belong to $K$, which is impossible.
The image of the generator $t_{21}^{(s+1)}$ under an
automorphism of the form \eqref{muphi}
is given by
$t_{21}^{(s+1)}+\vp_1\ts t_{21}^{(s)} +\cdots+\vp_{s}\ts t_{21}^{(1)}$.
Thus, the Verma module
$M(\la_1(u),\la_2(u))$ over $\Y(\gl_2)$ contains a nonzero
vector of the form
\ben
\ze=c_1\ts t_{21}^{(1)}\ts 1_{\la}+\cdots+c_{s}\ts t_{21}^{(s)}\ts 1_{\la}
+t_{21}^{(s+1)}\ts 1_{\la},
\qquad c_i\in\CC,
\een
such that $t_{12}^{(r)}\ts \ze=0$
for all $r\geqslant 1$.
The proof is completed by
writing this vector in terms of the elements $f_1^{(r)}$
with the use of the isomorphism \eqref{isom} and noting that
$t_{12}(u)\ts\ze=0$ implies $e(u)\ts\ze=0$.
\epf

\bre\label{rem:singular}
Nonzero vectors $\ze\in M(\mu(u))$ satisfying
$e(u)\ts\ze=0$ are analogous to the singular vectors
in the Verma modules over semisimple Lie algebras.
Note, however, that in contrast with the case of Lie algebras,
if $M(\mu(u))$ is reducible then, due to Corollary~\ref{cor:weaksing},
the subspace spanned by the vectors $\ze$ is infinite-dimensional.
Moreover, $\ze$ does not have to be an eigenvector
for the action of $h(u)$. Therefore, $\ze$
does not have to generate a submodule in $M(\mu(u))$ isomorphic
to a Verma module over the Yangian $\Y(\sll_2)$.
\qed
\ere

We now turn to the case of the Yangian $\Y(\agot)$ for
an arbitrary simple Lie algebra $\agot$ over $\CC$.
We start by proving the ``if" part of Theorem~\ref{thm:redu}.

\bpr\label{prop:ifpart}
Let $\mu(u)=(\mu_1(u),\dots,\mu_n(u))$ be the highest weight
of the Verma module $M(\mu(u))$ over the Yangian $\Y(\agot)$
such that for some $i\in\{1,\dots,n\}$ the series $\mu_i(u)$
is the Laurent expansion at $u=\infty$ of a rational function in $u$.
Then the module $M(\mu(u))$ is reducible.
\epr

\bpf
Consider the subalgebra of $\Y(\agot)$ generated by
the elements $e_i^{(r)}, h_i^{(r)}, f_i^{(r)}$ with $r=0,1,2,\dots$.
By the Yangian defining relations and the Poincar\'e--Birkhoff--Witt theorem
for $\Y(\agot)$, this subalgebra is isomorphic to the Yangian
$\Y(\sll_2)$. Moreover, the $\Y(\sll_2)$-span of the highest vector
$1_{\mu}$ of $M(\mu(u))$ is isomorphic to the Verma module $M(\mu_i(u))$ over
$\Y(\sll_2)$. By Corollaries~\ref{cor:sltwo} and
\ref{cor:weaksing}, this module contains a nonzero vector
of the form
\ben
\ze=c_0\ts f_{i}^{(0)}\tss 1_{\mu}+\cdots+c_p\ts f_{i}^{(p)}\tss 1_{\mu},
\een
such that $e_i(u)\ts\ze=0$. By the relations \eqref{ef},
we have $e_j(u)\ts\ze=0$ for all $j=1,\dots,n$. The
Poincar\'e--Birkhoff--Witt theorem
for $\Y(\agot)$ implies that the submodule $\Y(\agot)\ts\ze$ of
$M(\mu(u))$ does not contain the highest vector $1_{\mu}$.
Thus, $M(\mu(u))$ is reducible.
\epf

The definition of the Yangian $\Y(\agot)$ can be extended
to the case of semisimple Lie algebras $\agot$; see Section~\ref{sec:def}.
If $\agot=\agot_1\oplus\cdots\oplus\agot_p$ is a decomposition
of $\agot$ into the direct sum of simple ideals, then extending
the Poincar\'e--Birkhoff--Witt theorem to the Yangian $\Y(\agot)$
as in \cite{l:pb}, one can show that $\Y(\agot)$ is isomorphic
to the tensor product of the Yangians $\Y(\agot_1),\dots,\Y(\agot_p)$.

The proof of Theorem~\ref{thm:redu} is completed by the
following proposition which is the ``only if" part of the theorem.
We assume here that $\agot$ is a semisimple Lie algebra.

\bpr\label{prop:onlyifpart}
Suppose that the Verma module $M(\mu(u))$ over the Yangian $\Y(\agot)$
is reducible. Then there exists an index $i\in\{1,\dots,n\}$
such that the $i$-th component $\mu_i(u)$ of the highest weight
is the Laurent expansion at $u=\infty$ of a rational function in $u$.
\epr

\bpf
Consider the weight subspace decomposition \eqref{weidec}
of $M(\mu(u))$ with
respect to the Cartan subalgebra $\h$ of $\agot$.
Let us equip the set of weights with the standard partial
ordering: $\eta$ precedes $\eta^{\ts\prime}$ if $\eta^{\ts\prime}-\eta$
is a linear combination of the positive roots with nonnegative integer
coefficients. Any nontrivial submodule $K$ of $M(\mu(u))$ inherits
the weight space decomposition.
Taking a nonzero element of the weight subspace of $K$ with
a maximal weight $\eta$ we conclude that $M(\mu(u))$ must contain
a nonzero weight vector $\ze$ such that
\ben
e_i(u)\ts\ze=0,\qquad i=1,\dots,n.
\een
Let us write the weight $\eta$
in the form \eqref{weta} and suppose
that $k_1=\dots=k_{i-1}=0$
and $k_i>0$ for some index $i$. We shall be proving that
the $i$-th component $\mu_i(u)$ of the highest weight
has the desired property.
We may assume without loss
of generality that $i=1$. Indeed, otherwise we can consider the Yangian
for the semisimple Lie algebra $\agot'$ associated with the Cartan matrix
obtained from $(a_{ij})$ by deleting the first $i-1$ rows and columns.
The vector $\ze$ may then be regarded as an element of the Verma module
$M(\mu'(u))$ over $\Y(\agot')$ with the highest weight
$\mu'(u)=(\mu_i(u),\dots,\mu_n(u))$.

Let us fix a total ordering $\prec$ on the set of positive roots
$\Delta^+$ such that $\al_i$ precedes $\al_j$ if $i>j$ and any
composite root precedes any simple root.
Now consider the ordering on the set
$
\{f_{\al}^{(r)}\ |\ \al\in\Delta^+,\ r \geqslant 0\}
$
defined by the rule: $f_{\al}^{(r)}$ precedes $f_{\be}^{(s)}$
if $\al\prec\be$, or $\al=\be$ and $r<s$.
Let us write $\ze$ as a linear combination of the basis vectors
\eqref{ordmon} with the chosen ordering.
We claim that none of the generators $f_{\al}^{(r)}$
corresponding to a composite root $\al$ can occur in the expansion of
$\ze$. In other words, our claim is that $\ze$ is a linear combination
of monomials of the form
\beql{monsr}
f_{\al_n}^{(p_{1_{}})}\cdots f_{\al_n}^{(p_{k_n})}\cdots
f_{\al_1}^{(q^{}_1)}\cdots f_{\al_1}^{(q^{}_{k_1})}\tss 1_{\mu}.
\eeq
Write $\ze=\ze_{\ts 0}+\cdots+\ze_{\ts m}$, where $\ze_{\ts k}$
is a linear combination
of basis monomials \eqref{ordmon} of degree $k$. We shall prove
by a reverse induction on $k$ that $\ze_{\ts k}$ is a linear combination
of monomials of the form \eqref{monsr}.
Let $\be$ be a positive root such that $\ga:=\be-\al_1$ is a root.
Suppose that $s$ is
the maximum nonnegative integer such that
the element $f_{\be}^{(s)}$ occurs in the expansion
of $\ze_{\ts k}$. By the construction of the elements $f_{\al}^{(r)}$,
we have
\ben
[e_{\al_1}^{(r)},f_{\be}^{(s)}]\equiv c\ts f_{\ga}^{(r+s)},
\een
modulo elements of degree smaller than $r+s$, where
$c$ is a nonzero constant. Therefore, the expansion of $e_{\al_1}^{(r)}\ts\ze$
into a linear combination of basis monomials will contain
a monomial of degree $r+k$ where $f_{\ga}^{(r+s)}$ occurs as a factor.
By the induction hypothesis, the components $\ze_{\ts k+1},\dots,\ze_{\ts m}$
of $\ze$ are linear combinations of monomials of the form \eqref{monsr}.
Hence, for a sufficiently large $r$
the expansion of $e_{\al_1}^{(r)}\ts\ze_j$ with $j=k+1,\dots,m$ will not contain
$f_{\ga}^{(r+s)}$.
Since $e_{\al_1}^{(r)}\ts\ze=0$ for all $r\geqslant 0$, we come to
a contradiction, as the monomial containing $f_{\ga}^{(r+s)}$
with a sufficiently large $r$
will occur in the expansion of $e_{\al_1}^{(r)}\ts\ze$ with a nonzero
coefficient. So,
if $\beta$ is a positive root such that $\be-\al_1$ is a root then
none of the elements $f_{\be}^{(s)}$ can occur in the expansion
of $\ze_{\ts k}$.

Similarly, if $\be$ is a composite positive root
such that $\be-\al_1$ is not a root, but $\be-\al_2$ is a root,
then we use the relations $e_{\al_2}^{(r)}\ts\ze=0$ to show that
$f_{\be}^{(s)}$ cannot occur in the expansion of $\ze_{\ts k}$ either.
Continuing in the same manner, we conclude that $\ze_{\ts k}$
may only contain the generators $f_{\be}^{(s)}$ where $\be$
is a simple root thus proving the claim.

Thus, we may write the vector $\ze$ as the sum of vectors
of the form
\ben
f_{\al_n}^{(p_{1_{}})}\cdots f_{\al_n}^{(p_{k_n})}\cdots
f_{\al_2}^{(s^{}_1)}\cdots f_{\al_2}^{(s^{}_{k_2})}\tss \ze',
\een
where each vector $\ze'$ is a linear combination
of monomials
$
f_{\al_1}^{(q^{}_1)}\cdots f_{\al_1}^{(q^{}_{k_1})}\tss 1_{\mu}
$
and at least one of the vectors is nonzero.
However, the relations $e_{\al_1}^{(r)}\ts\ze=0$
imply $e_{\al_1}^{(r)}\ts\ze'=0$ for each vector $\ze'$.
Recall that
the elements $e_{\al_1}^{(r)}$, $f_{\al_1}^{(r)}$ and $h_1^{(r)}$
generate a subalgebra of $\Y(\agot)$ isomorphic to the Yangian $\Y(\sll_2)$.
So, the Verma module $M(\mu_1(u))$ over $\Y(\sll_2)$
is reducible. By Corollary~\ref{cor:weaksing}, the series $\mu_1(u)$
is the Laurent expansion at $u=\infty$ of a rational function in $u$.
\epf

\section{Weight subspaces of the irreducible quotient}\label{sec:fin}
\setcounter{equation}{0}

We now prove Theorem~\ref{thm:fin}. The following is the ``only if"
part of the theorem.

\bpr\label{prop:onlyiffin}
Suppose that
all weight subspaces of the $\agot$-module $L(\mu(u))$
are finite-dimensional. Then each component $\mu_i(u)$
of the highest weight
is the Laurent expansion at $u=\infty$ of a rational function in $u$.
\epr

\bpf
Given index $i\in\{1,\dots,n\}$ identify the subalgebra
of $\Y(\agot)$ generated by the elements $e_i^{(r)}$, $f_i^{(r)}$
and $h_i^{(r)}$ for $r\geqslant 0$ with the Yangian $\Y(\sll_2)$.
The $\Y(\sll_2)$-span of the highest vector $1_{\mu}$ in $L(\mu(u))$
is isomorphic to a quotient $L$ of the Verma module $M(\mu_i(u))$ over $\Y(\sll_2)$.
By our assumptions, all weight subspaces of $L$, regarded as an $\sll_2$-module,
are finite-dimensional. Therefore, the module $M(\mu_i(u))$ is reducible.
The proof is completed by the application of Corollary~\ref{cor:sltwo}.
\epf

The following proposition completes the proof of Theorem~\ref{thm:fin}.

\bpr\label{prop:iffin}
Suppose that each component $\mu_i(u)$
of the highest weight of the $\Y(\agot)$-module $L(\mu(u))$
is the Laurent expansion at $u=\infty$ of a rational function in $u$.
Then all weight subspaces of the $\agot$-module $L(\mu(u))$
are finite-dimensional.
\epr

\bpf
Each component $\mu_i(u)$ is the Laurent expansion at $u=\infty$
of a rational function $P_i(u)/Q_i(u)$, where $P_i(u)$ and $Q_i(u)$
are monic polynomials in $u$ of the same degree. Let $p_i$ be the degree
of these polynomials. For any positive root $\al$ denote by
$[\al:\al_i]$ the multiplicity of the simple root
$\al_i$ in $\al$.
We shall be proving that the monomials of the form
\eqref{ordmon} with the condition
\beql{degrj}
r_j<\sum_{i=1}^n [\al^{(j)}:\al_i]\ts p_i,\qquad j=1,\dots,k,
\eeq
span the module $L(\mu(u))$. This clearly implies that each weight subspace
$L(\mu(u))_{\eta}$ is finite-dimensional.

It suffices to show that each monomial \eqref{ordmon} can be
written in $L(\mu(u))$ as a linear combination of those
monomials satisfying \eqref{degrj}. We argue by induction on
the degree $r$ of the monomial. With this degree fixed, we also use
induction on the length of the monomial.
As the induction base, consider
a monomial $f_{\be}^{(r)}\ts 1_{\mu}$ of degree $r$ and length $1$,
where $\be$ is a positive root. There is nothing to prove unless
$r$ satisfies
\beql{rcond}
r\geqslant \sum_{i=1}^n [\be:\al_i]\ts p_i.
\eeq
In this case, by definition of the elements $f_{\be}^{(r)}$, we can write,
modulo terms of smaller degree,
\beql{fbe}
f_{\be}^{(r)}\equiv [f_{i_1}^{(r_{i_1})},[f_{i_2}^{(r_{i_2})},\dots,
[f_{i_{l-1}}^{(r_{i_{l-1}})},f_{i_l}^{(r_{i_l})}]\dots ]],
\eeq
for a partition $r=r_{i_1}+\cdots+r_{i_l}$
which due to \eqref{rcond} can be chosen in such
a way that $r_{i_a}\geqslant p_{i_a}$
for all $a$. However, for any index $i\in\{1,\dots,n\}$
and any $r_i\geqslant p_i$ the monomial $f_{i}^{(r_i)}\ts 1_{\mu}$
is a linear combination of monomials of smaller degree. Indeed,
this follows from Corollary~\ref{cor:weaksing} by
identifying
the subalgebra
of $\Y(\agot)$ generated by the elements $e_i^{(r)}$, $f_i^{(r)}$
and $h_i^{(r)}$ for $r\geqslant 0$ with the Yangian $\Y(\sll_2)$.
Therefore, expanding the commutators in \eqref{fbe} we conclude
that $f_{\be}^{(r)}\ts 1_{\mu}$ is a linear combination
of monomials of smaller degree.

Now consider a monomial of the form \eqref{ordmon}
of degree $r$, where the condition \eqref{degrj} is violated
for some $j$. It will suffice to assume that $j=1$.
Recall that modulo terms of smaller degree,
\ben
[f_{\al}^{(r)}, f_{\be}^{(s)}]\equiv c\ts f_{\al+\be}^{(r+s)},\qquad c\in\CC,
\een
if $\al+\be$ is a root; otherwise the commutator is zero.
Using this relation,
we shall move the factor $f_{\al^{(1)}}^{(r_1)}$ in \eqref{ordmon}
to the right so that the monomial will be equal,
modulo terms of smaller degree or smaller length, to the monomial
\ben
f_{\al^{(2)}}^{(r_2)}\cdots f_{\al^{(k)}}^{(r_k)}\ts
f_{\al^{(1)}}^{(r_1)}\ts 1_{\mu}.
\een
However, by the induction base, $f_{\al^{(1)}}^{(r_1)}\ts 1_{\mu}$
is zero, modulo terms of smaller degree.
\epf

\bre\label{rem:character}
In the case where all weight subspaces of $L(\mu(u))$ are finite-dimensional
we can consider the character of the $\agot$-module $L(\mu(u))$,
\ben
\text{ch\ts} L(\mu(u))=\sum_{\eta} \dim L(\mu(u))_{\eta}\cdot e^{\eta},
\een
where $e^{\eta}$ is a formal exponential; see e.g. Dixmier~\cite[Chapter~7]{d:ae}.
It would be interesting to find a formula for this character.
Character formulas for finite-dimensional modules $L(\mu(u))$ are given
e.g. in \cite{a:df} and \cite{v:aq}. In the particular case $\agot=\sll_2$,
a character formula
for $L(\mu(u))$ can be deduced from the tensor product decomposition
for $L(\mu(u))$; see \cite{t:im} and \cite[Proposition~3.6]{m:fd}.
More precisely, write the rational function $\mu(u)$ as
\ben
\mu(u)=\frac{(u+\al_1)\cdots (u+\al_k)}{(u+\be_1)\cdots (u+\be_k)},\qquad
\al_i,\be_i\in\CC.
\een
Renumbering the $\al_i$ and $\be_i$ if necessary, we may assume
without loss of generality that
for every $i=1,\dots, k-1$ the following condition holds:
if the multiset
\ben
\{\al_p-\be_q\ |\ i\leqslant p,q \leqslant k\}
\een
contains nonnegative integers, then $\al_i-\be_i$
is minimal amongst them. Let $l$ be the number of the indices $i$ such that
$\al_i-\be_i$ is a nonnegative integer. Then
\ben
\text{ch\ts} L(\mu(u))=
\prod_{i=1}^l
\frac{x^{\ts\al_i-\be_i+1}-x^{-\al_i+\be_i-1}}{x-x^{-1}}
\prod_{i=l+1}^k
\frac{x^{\ts\al_i-\be_i+1}}{x-x^{-1}},
\een
where $x=e^{\omega}$ and $\omega$ is the fundamental
weight for $\sll_2$.
\ere

\bre\label{rem:qanalog}
The results of this paper can be extended to obtain
analogs of Theorems~\ref{thm:redu} and \ref{thm:fin}
for the quantum affine algebras $\U_q(\wh\agot)$;
see e.g. \cite{d:ha, d:nr}
for their definition.
\ere

\section*{Acknowledgments}

The first author gratefully acknowledges the support from
the Natural Sciences and Engineering Research Council of Canada.
The second author is
supported by the
CNPq grant (processo 307812/2004-9) and the Fapesp
grant (processo 04/06258-0).
The third
author is grateful to the Fapesp for the financial support
(processo 04/10038-6). He would like to thank the University of S\~ao
Paulo and the Carleton University
for the warm hospitality during his visits.


\begin{thebibliography}{99}

\bibitem{a:df}
{T. Arakawa}, {\it Drinfeld functor and finite-dimensional
representations of Yangian\/},
Comm. Math. Phys. {\bf 205} (1999), 1--18.

\bibitem{bz:wm}
{Y. Billig and K. Zhao},
{\it Weight modules over exp-polynomial Lie algebras},
J. Pure and Appl. Algebra {\bf 191} (2004), 23--42.

\bibitem{cp:gq}
{V. Chari and A. Pressley},
{\it A guide to quantum groups},
Cambridge University Press, 1994.

\bibitem{d:ae}
{J. Dixmier},
{\it Alg\`ebres Enveloppantes},
{Gauthier-Villars, Paris},
1974.

\bibitem{d:ha}
{V. G. Drinfeld},
{\it Hopf algebras and the
quantum Yang--Baxter equation},
{Soviet Math. Dokl.} {\bf
32} (1985), 254--258.

\bibitem{d:nr}
{V. G. Drinfeld},
{\it A new realization of
Yangians and quantized affine algebras}, {Soviet Math. Dokl.}
{\bf 36} (1988), 212--216.

\bibitem{l:pb}
{S. Z. Levendorski\u\i},
{\it On PBW bases for Yangians},
{Lett. Math. Phys.}
{\bf 27}
(1993),
37--42.

\bibitem{m:fd}
{A. I. Molev},
{\it Finite-dimensional irreducible representations of twisted
Yangians},
{J. Math. Phys.} {\bf 39} (1998), 5559--5600.

\bibitem{mno:yc}
{A. Molev, M. Nazarov and G. Olshanski},
{\it Yangians and classical Lie algebras},
Russian Math. Surveys
{\bf 51}:2
(1996),
205--282.

\bibitem{t:im}
{V. O. Tarasov},
{\it Irreducible monodromy matrices for the $R$-matrix of the
$XXZ$-model and lattice local quantum Hamiltonians}, {Theor. Math. Phys.}
{\bf 63} (1985),
440--454.

\bibitem{v:aq}
{E. Vasserot},
{\it Affine quantum groups and equivariant $K$-theory},
Transform. Groups {\bf 3} (1998),
269--299.


\end{thebibliography}
\end{document}